\documentclass[11pt]{article}
\usepackage{amsmath,amsthm,amsfonts,latexsym,amssymb,enumerate,color}
\usepackage{graphicx}
\usepackage{mathrsfs}

\def\CC{\mathbb C}
\def\DD{\mathbb D}
\def\EE{\mathbb E}
\def\GE{\mathcal E}
\def\GF{\mathcal F}
\def\G{\mathcal G}
\def\GH{\mathcal H}
\def\GK{\mathcal K}
\def\GL{\mathcal L}
\def\GM{\mathcal M}
\def\GN{\mathcal N}
\def\NN{\mathbb N}
\def\GP{\mathcal P}

\def\GS{\mathcal S}
\def\TT{\mathbb T}

\def\ker{\mathop{\rm ker}\nolimits}

\def\spam{\mathop{\rm span}\nolimits}  

\newtheorem{thm}{Theorem}[section]
\newtheorem{prop}[thm]{Proposition}

\newtheorem{cor}[thm]{Corollary}
\newtheorem{rem}[thm]{Remark}

\numberwithin{equation}{section}

\def\beginpf{\begin{proof}}
\def\endpf{\end{proof}}
\def\beq{\begin{equation}}
\def\eeq{\end{equation}}

\def\ol{\overline}

\def\sab{S_{a,b}}

\def\k{{\rm ker}}
\def\kerab{\k_{a,b}}

\def\H2b{\overline{H^2_0}}
\def\Kmin{K_{\rm min}}

\begin{document}

\title{Paired kernels and their applications}

\author{M.~Cristina C\^amara,\thanks{
Center for Mathematical Analysis, Geometry and Dynamical Systems,
Instituto Superior T\'ecnico, Universidade de Lisboa, 
Av. Rovisco Pais, 1049-001 Lisboa, Portugal.
 \tt ccamara@math.ist.utl.pt} 
 \and  Jonathan R.~Partington\thanks{School of Mathematics, University of Leeds, Leeds LS2~9JT, U.K. {\tt j.r.partington@leeds.ac.uk}
}
\  }

\date{}

\maketitle

\begin{abstract}
This paper considers paired operators in the context of the Lebesgue
Hilbert space on the unit circle and its subspace, the Hardy space $H^2$.
The kernels of such operators, together with
their analytic projections, which are generalizations of Toeplitz kernels, are studied. 
Results on
near-invariance properties, representations, and inclusion relations for these kernels
are obtained. The existence of a minimal Toeplitz kernel containing any projected paired kernel
and, more generally, any nearly $S^*$-invariant subspace of $H^2$, is derived.
The results are applied to describing the kernels of finite-rank asymmetric
truncated Toeplitz operators.
\end{abstract}


\footnote{This work was partially
supported by FCT/Portugal through CAMGSD, IST-ID, projects UIDB/04459/2020 and UIDP/04459/2020.}

\section{Introduction}

Let $X$ be a Banach space, $P \in \GL(X)$ a projection, and $Q=I-P$ the complementary projection. An operator of the form $AP+BQ$ or $PA+QB$ with $A,B \in \GL(X)$ is called
a {\em paired operator\/}
\cite{CG,MP,prossdorf}. In this paper we consider the case when $A=M_a$ and $B=M_b$ are
multiplication operators on $L^2:=L^2(\TT)$, where $\TT$ denotes the unit circle, with $a, b\in L^\infty:=L^\infty(\TT)$, and we denote by $S_{a,b}$ and $\Sigma_{a,b}$ the operators defined on $L^2$ by
\beq\label{S}
S_{a,b}f= aP^+f+bP^-f \,\,, \qquad \Sigma_{a,b} f = P^+(af)+P^-(bf).
\eeq

Paired operators first appeared in the context of the theory of
singular integral equations 
\cite{widom,shinbrot}. Consider the canonical example of a singular integral operator on $L^2$,
\beq
(Lf)(t) = A(t)f(t) + B(t) (S_\TT f)(t),
\eeq
with $A,B \in L^\infty$ and
\beq
 (S_\TT f)(t) = \frac{1}{\pi i} \hbox{PV} \int_\TT \frac{f(y)}{y-t} \, dy, \qquad t \in \TT.
 \eeq

 It is well known that, denoting by $P^\pm$ the orthogonal
 projections from $L^2$ onto the Hardy spaces $H^2_+:=H^2(\DD)$ and
 $H^2_- = \H2b=(H^2_+)^\perp$, respectively, 
 identified with closed subspaces of $L^2$,
 we have $I=P^+ + P^-$ and $S_\TT=P^+-P^-$, so that
 the operator $L$ can be expressed as
 \beq
 L= (A+B)P^+ + (A-B)P^-. 
 \eeq
We may write this as  $S_{a,b} $
with $ a=A+B\,,\,b=A-B$,
 while its adjoint is a paired operator of the second type in \eqref {S}, 
 \beq
S^*_{a,b}f=\Sigma_{\bar a,\bar b} .
 \eeq

Paired operators are also closely related to
Toeplitz operators of the form
\beq
T_a = P^+ a P^+ _{|H^2},
\eeq
where $a \in L^\infty$ is called the {\em symbol} of the operator. If we represent $S_{a,b}$ in the form
\beq\label{dil}
\begin{pmatrix}
P^+aP^+  & P^+ b P^- \\
P^- aP^+ & P^-b P^- 
\end{pmatrix}: H^2_+ \oplus H^2_- \to H^2_+ \oplus H^2_-,
\eeq
we see that paired operators are dilations of Toeplitz operators. If $b \in \G L^\infty$ 
 (that is, invertible in $L^\infty$)
we can write
\beq
S_{a,b} = aP^+ + bP^- = b\left( \frac{a}{b}P^+ + P^- \right),
\eeq
which is equivalent after extension \cite{BTsk} to the Toeplitz operator $T_{a/b}$
\cite{C}.

 However, paired operators and spaces that turn out to be kernels of paired operators, called {\em paired kernels}, also
appear in different guises, for instance in the
study of dual truncated Toeplitz operators \cite{CKP}, in the description of scalar-type block Toeplitz
kernels \cite{CP20}, in the characterization of the ranges of
finite-rank truncated Toeplitz operators \cite{CC19}, and in the study of nearly invariant subspaces for shift semigroups \cite{LP22}.

 We shall consider mainly paired operators of the form $S_{a,b}$, where the
pair $(a,b)$ is called a {\em symbol pair}. These operators
have been considered mostly under particular conditions, such as invertibility in $L^\infty$, for $a,b$ or $a/b$ \cite{widom,shinbrot}; see also \cite{MP} and references therein.
The operator $S_{a,b}$ is said to be of {\em normal type} if $a,b \in \G L^\infty$.
This is by far the most studied case, but some particular types of
non-normal paired operators have also been considered  \cite{MP,prossdorf,CMP21,castrigiano}.
We shall assume throughout the paper, more generally, that $a,b \in L^\infty$ and
the pair $(a,b)$ is {\em nondegenerate}; that is, $a,b$ and $a-b$ are nonzero a.e.\ on $\TT$.
We also use the notation $P^\pm \phi = \phi_\pm$ for $\phi \in L^2$.

Here we shall study in particular the properties of kernels of
paired operators, called {\em paired kernels}, and their projections into
$H^2_+$ and $H^2_-$, called {\em projected paired kernels}. Denoting
 \beq
 \kerab = \ker \sab, \qquad \k^\pm_{a,b}= P^\pm \kerab,
 \eeq
 we have that $\ker T_{a/b} = P^+ \ker S_{a,b} =: \k^+_{a,b}$ for the possibly unbounded Toeplitz operator $T_{a/b}$. If $a/b\in L^\infty$ then $ \k^+_{a,b}$ is a {\em Toeplitz kernel}, i.e., the kernel of a bounded Toeplitz operator; otherwise it may not be a closed subspace of $H^2_+:=P^+L^2$. However, we can define a one-to-one correspondence (see Section 2) between  $ \k^+_{a,b}$ and $\kerab$, where the latter is closed because it is the kernel of a bounded operator. So one can also see paired kernels as being the natural closed space generalizations of Toeplitz kernels, allowing us to study the
kernels of unbounded Toeplitz operators
of the form $T_{a/b}$ with $a/b \not\in L^\infty$ in terms of the bounded operators
$S_{a,b}$ and $P^+$ on $L^2$.
It is thus natural to ask whether some known properties
of Toeplitz kernels with bounded symbols
can be extended or related with corresponding properties of paired kernels
or their projections on $H^2_+$. Alternatively one can look at this as studying the question of how certain properties of a Toeplitz operator extend to a dilation of the form \eqref{dil}.

In the following sections we study several properties of paired kernels which extend or are in contrast with various known properties of Toeplitz kernels. To compare the case where $a/b\in L^\infty$ (and $\k^+_{a,b}$ is a Toeplitz kernel) with that where $a/b\notin L^\infty$ and to illustrate some natural questions arising in the latter case, we start by considering in Section \ref{sec:3} an example which appears in the study of nearly invariant subspaces for shift semigroups \cite{LP22}. This leads to the study of near invariance properties of paired kernels (Section \ref{sec:n4}); to the question of existence of a minimal Toeplitz kernel containing any given $\k^+_{a,b}$ 
and, more generally, any space of the form $b \ker T_a$, including the
closed nearly $S^*$-invariant subspaces of $H^2_+$
(Section \ref{sec:n5}); to investigating the relations between paired kernels for operators with connected symbol pairs (Section \ref{sec:6}). In Section \ref{sec:7} the results are applied to study and describe the kernels of finite rank asymmetric truncated Toeplitz operators \cite{CP17}, showing in particular that, surprisingly, they do not depend on the range space if the latter is ``large" enough.


\section{Projected paired kernels}

Recall that we denote the kernel of a paired operator, which we call a {\em paired kernel}, by
 \beq
 \kerab = \ker \sab
 \eeq
and we call
 \beq\label{eq:3.5}
 \k^\pm_{a,b}= P^\pm \kerab = P^\pm \ker \sab
 \eeq
a {\em projected paired kernel}.

Each function $f \in L^2$ belongs to one and only one
paired kernel \cite[Thm. 4.6]{CGP23}. On the other hand,
for any $\phi_+ \in H^2_+$,
\[
\phi_+ \in \kerab^+ \iff \phi_+ \in \ker T_{a/b},
\]
where $T_{a/b}$ is the possibly unbounded Toeplitz operator with symbol $a/b$ defined on the domain
\[
D_{a/b} = \{ \phi_+ \in H^2_+: \frac{a}{b}\phi_+ \in L^2 \}.
\]
Analogously, for the dual Toeplitz operator \cite{CKP}
\[
\check T_{b/a} : \check D_{b/a}= \left\{ \phi_- \in H^2_-: \frac{b}{a} \phi_- \in L^2 \right\} \to H^2_-,
\]
defined by 
\[
\check T_{b/a} \phi_- = P^- \frac{b}{a} \phi_-,
\]
we have that $\kerab^- = \ker \check T_{b/a}$.

If $a/b \in L^\infty$ then $D_{a/b} = H^2_+$ and we say that $\kerab^+$ is a {\em Toeplitz kernel}, i.e., the kernel of a bounded Toeplitz operator.

Paired kernels and their projections into $H^2_\pm$ can be related as follows.

\begin{prop}\label{nprop:4.2}
The operators
\beq
\GP^+ :\kerab \to \kerab^+, \qquad \GP^+ \phi =  \frac{b}{b-a} \phi
\eeq
and
\beq
\GM_{a/b}: \kerab^+ \to \kerab^-, \qquad \GM_{a/b} \phi_+ = \frac{a}{b}\phi_+
\eeq
are well-defined and bijective with inverses
\beq
({\GP^+})^{-1}: \kerab^+ \to \kerab\,\,, \qquad ({\GP^+})^{-1}\phi_+ = (1-\frac{a}{b})\phi_+
\eeq
and
\beq
\GM^{-1}_{a/b}: \kerab^- \to \kerab^+\,\,, \qquad \GM^{-1}_{a/b} \phi_-= \frac{b}{a} \phi_-,
\eeq
and we have $\GP^+ \phi = P^+ \phi$, and $\GM_{a/b}(\GP^+ \phi)=P^- \phi$
for $\phi \in \kerab$.
\end{prop}

\begin{cor}\label{cor:2.2}
$\dim \kerab < \infty \iff \dim \kerab^+ < \infty \iff \dim \kerab^- < \infty$
and, if these dimensions are finite, then they are equal.
\end{cor}

The following is also an immediate consequence of Proposition \ref{nprop:4.2},
noting that if $\phi_+=0$ (or similarly for
 $\phi_-$) on a set of positive measure then, by the Luzin--Privalov theorem, it is $0$ a.e. on $\TT$.

\begin{cor}\label{cor:2.3bis}
If $\phi \in \kerab$ then $\phi=0 \iff \phi_+=0 \iff \phi_-=0$.
\end{cor}

Note that, although $\GP^+$ is bounded and bijective, its inverse is not necessarily bounded when $b \not\in \G L^\infty$, since $\kerab^+$ may not be closed in $H^2_+$, as shown in 
an example in the next section.

Just as we can relate the kernels of Toeplitz operators with those of dual Toeplitz operators, we can also reduce the study of $\kerab^-$ to that of $\k^+_{\bar b,\bar a}$, as the next proposition shows.

\begin{prop}\label{nprop:2.4}
$\kerab^- = \bar z \overline{\k^+_{\bar b,\bar a}}$.
\end{prop}

\beginpf For $\phi_-\in H^2_-$ we have
\begin{eqnarray*}
\phi_- \in \kerab^-  & \iff & a\phi_++b\phi_- = 0 \hbox{ for some } \phi_+ \in H^2 \\
& \iff & \bar a(\bar z \overline{\phi_+}) + \bar b(\bar z \overline{\phi_-})=0 \hbox { for some } \phi_+ \in  H^2 \\
& \iff & \bar b(\bar z \overline{\phi_-}) + \bar a \psi_- =0 \hbox { for some } \psi_- \in \H2b \\
& \iff & \bar z \overline{\phi_-} \in \k^+_{\bar b,\bar a} \iff \phi_- \in \bar z \overline{\k^+_{\bar b, \bar a}}.
\end{eqnarray*}
\endpf

From now on we shall mainly focus on the properties of the projected
paired kernels $\kerab^+$.
To compare the case where $b \in \G L^\infty$ with that where $b \not \in \G L^\infty$, we start by considering a particular example
of the latter, which is used in \cite{LP22}.


\section{An example and questions it raises}\label{sec:3}

Let $\theta$ be the singular inner function $\theta(z)=\exp \dfrac{z-1}{z+1}$ for $z \in \TT$
and let $a=\bar\theta$, $b(z)=z+1$. The kernel of $S_{\bar\theta,z+1}$ is described by
\beq
\bar\theta \phi_+ + (z+1)\phi_-=0, \qquad \hbox{i.e.,} \quad\quad \frac{\phi_+}{z+1}=-\theta \phi_-.
\eeq
Since the left-hand side of this last equation represents a function in the
Smirnov class $\GN_+$ and the right-hand side represents
a function in $L^2$, we have that both belong to $H^2_+$ and thus, from
\beq
\bar\theta \frac{\phi_+}{z+1} = -\phi_-,
\eeq
we conclude that $\phi_+/(z+1) \in K_\theta$, where $K_\theta$
denotes the model space $ H^2_+ \ominus \theta H^2_+ = \ker T_{\bar\theta}$.
It follows that $\k^+_{\bar\theta,z+1}\subseteq (z+1)K_\theta$ and the converse inclusion is easily seen to be true, so
\beq
\k^+_{\bar\theta,z+1}=(z+1)K_\theta.
\eeq
Clearly $\ker_{\bar\theta,z+1}=\left(1-\frac{\bar\theta}{z+1} \right) \k^+_{\bar\theta,z+1}$
is a closed subspace of $L^2$, since $S_{\bar\theta,z+1}$ is bounded, but
\beq
\GP^+ \k_{\bar\theta,z+1} = \k^+_{\bar\theta,z+1}=(z+1)K_\theta
\eeq
is not closed \cite[Prop 3.2]{LP22}, so $\GP^+$ does not have a bounded inverse.

Note that, while $\k^+_{\bar\theta,z+1}$ is not closed, it is nevertheless contained
in a (closed) minimal Toeplitz kernel, by which we mean a Toeplitz kernel that contains 
$\k^+_{\bar\theta,z+1}$ and it itself contained in any other Toeplitz kernel containing
$\k^+_{\bar\theta,z+1}$. Indeed, on the one hand,
\beq\label{neq:3.5}
(z+1)K_\theta = \k^+_{\bar\theta,z+1} \subsetneq K_{z\theta}=\ker T_{\bar z \bar\theta}.
\eeq
On the other hand it was shown in \cite {CP14} that for any nonzero $\phi_+\in H^2_+$ there exists a minimal Toeplitz kernel to which $\phi_+$ belongs. Now, $\ker T_{\bar z\bar\theta}$ is the minimal kernel
containing the function 
\[
f = (z+1) \frac{\theta-\theta(0)}{z} \in \k^+_{\bar\theta,z+1},
\]
since we have 
\[
\bar z\bar\theta f = \bar z \frac{z+1}{z} (1-\theta(0)\bar\theta) = \bar z \overline{(z+1)(1-\overline{\theta(0)}\theta)},
\]
where $(z+1)(1-\overline{\theta(0)}\theta)$ is an outer function in $H^2_+$
(see \cite[Thm, 2.2]{CP18a}). Thus $K_{z\theta  }$ is the minimal Toeplitz kernel containing
$\k^+_{\bar\theta,z+1}$.

These results naturally raise several questions,
especially when compared with some known results for Toeplitz kernels.\\

{\bf Question \ref{sec:3}.1.} 
Having shown that $\k^+_{\bar\theta,z+1}$ is not a Toeplitz kernel, one may ask
whether 
 there is any Toeplitz kernel contained in that space. 
The answer is negative, due to the
near invariance properties of Toeplitz kernels
\cite{CP14}; these imply
in particular that no Toeplitz kernel can be contained in $(z+1)H^2_+$, since Toeplitz kernels are
nearly $\dfrac{1}{z+1}$   invariant \cite{CP14}
and $\k^+_{\bar\theta,z+1}=(z+1)K_\theta$.
This equality also shows that, in contrast with Toeplitz kernels, projected paired kernels may not be
nearly $\dfrac{1}{z+1}$   invariant. 
But do other near invariance properties of Toeplitz kernels extend to
projected paired kernels? This is studied in Section \ref{sec:n4}, comparing the two cases
where $a/b \in L^\infty$ and $a/b \not\in L^\infty$.\\

{\bf Question \ref{sec:3}.2.} 
On the other hand, there exists a minimal Toeplitz kernel containing $\k^+_{\bar\theta,z+1}$,
which is $K_{\theta z}$. Is there a 
minimal Toeplitz kernel containing any given nontrivial 
$\kerab^+$?
We   answer   this question in the affirmative, and in a more general setting, in Section \ref{sec:n5}, by showing that
the closure of any projected paired kernel admits a representation of the form
$f \ker T_g$ with $f \in H^2_+$ and $g \in L^\infty$, and we discuss the
existene of such a representation for
projected paired kernels. Note that not only can every Toeplitz kernel be
writen as a product of the form $f \ker T_g$ \cite{hayashi}, but we also
have the same property for other projected paired kernels: for instance,
$\k^+_{\bar\theta,z+1}=(z+1)K_\theta$.\\

{\bf Question \ref{sec:3}.3.}
The relation \eqref{neq:3.5} can be rewritten as $\k^+_{\bar\theta,z+1} \subsetneq\k^+_{\bar\theta \bar z,1}$, raising the question of what may be
the inclusion relations between the two projected paired
kernels when we multiply the two elements of the symbol pair $(a,b)$ by certain functions.
For Toeplitz operators we have, for instance, that
\[
h_- \in \overline{H^\infty} \implies \ker T_g \subseteq \ker T_{h_- g},
\]
where the inclusion is strict if the inner factor of $\overline{h_-}$ is non-constant
\cite{CMP16} and an equality if $\overline{h_-}$ is outer in $H^\infty$; also,
\[
h_+ \in H^\infty \implies h_+ \ker T_{h_+ g} \subseteq \ker T_g,
\]
where the inclusion is strict if $h_+$ has a non-constant inner factor \cite{CMP16}
and an equality if $h_+$ is invertible in $H^\infty$.
We study how analogous inclusion relations can be established
for general projected paired kernels, and whether the inclusion is strict or an equality,
in Section \ref{sec:6}.

The results are applied to study kernels of (asymmetric) truncated Toeplitz operators in Section \ref{sec:7}.


\section{Near invariance properties} \label{sec:n4}

A subspace $\GS \subseteq H^2_+$ is said to be
{\em nearly $S^*$-invariant\/} if and only if
\beq \label{neq:5.1}
  f_+ \in \GS, f_+(0)=0 \implies S^*f_+ \in \GS.
  \eeq
 Here $S^*$ denotes the backward shift on $H^2_+$, i.e.,
  $S^*=T_{\bar z}$.  
  
  Noting that $f_+(0)=0$ means that $\bar z f_+ \in H^2_+$ and, in that case, $S^* f=\bar z f$, it is clear that
  \eqref{neq:5.1} is equivalent to
  \beq
  f_+ \in \GS, \bar z f_+ \in H^2_+ \implies \bar z f_+ \in \GS 
  \eeq
  and we say, equivalently, that $\GS$ is nearly $\bar z$-invariant \cite{CP14}. More generally, if $\eta$ is a complex-valued
  function defined a.e.\ on $\TT$ we say that $\GS \subseteq H^2_+$ is {\em nearly $\eta$-invariant\/} if and only if
  \beq
    f_+ \in \GS, \eta f_+ \in H^2_+ \implies \eta f_+ \in \GS.
    \eeq
    In this case if $\eta \in L^\infty$  we can also say that $\GS$ is nearly $T_\eta$-invariant. \\
    
    Toeplitz kernels are closed nearly $S^*$-invariant spaces. Furthermore, in \cite{CP14} a large class of
    functions $\eta$ was described for which all Toeplitz kernels are nearly
    $\eta$-invariant, and which includes all functions
    in $\ol{H^\infty}$ and all rational functions without poles in $\DD^e \cup \{\infty\}$,
    where $\DD^e = \{z \in \CC: |z|>1 \}$.
    
    In particular, Toeplitz kernels are nearly $\bar\theta$-invariant and nearly $\frac{1}{z-z_0}$-invariant, where $\theta$ is any inner function and $z_0 \in \TT\cup\DD$.
    As a consequence of this, we conclude that no Toeplitz kernel can be contained in $\theta H^2_+$ or in
    $(z-z_0)H^2_+$ for $z_0 \in \TT \cup \DD$.

It is clear from the examples of Section \ref{sec:3} that the latter property cannot be extended to
projected paired kernels in general. Other near invaiance properties of Toeplitz kernels, however, are
shared with
projected paired kernels.  

The following results are simple consequences of the definitions at the beginning
of this section. We assume that $\kerab^+ \ne \{0\}$.

\begin{prop}\label{prop:n4.1}
 $\kerab^+$ is nearly $\eta$-invariant for every $\eta \in \ol{H^\infty}$.
\end{prop}

\begin{cor}\label{ncor:5.3}
$\kerab^+$ is nearly $\bar\theta$-invariant for every inner function $\theta$,
and therefore there exists a function $\phi_+ \in \kerab^+$ such that
$\phi_+ \not\in \theta H^2_+$.
\end{cor}

\begin{cor}
$\k^+_{a,b}$ is nearly $R$-invariant for every rational function $R$ bounded at $\infty$ whose poles lie in $\DD$.
\end{cor}

Naturally, if $a/b \in L^\infty$, this property can be extended to every rational $R$ without
poles in $\DD^e \cup \{\infty\}$, since in that case $\k^+_{a,b}$ is a Toeplitz kernel.\\

It follows from Proposition \ref{prop:n4.1}, in particular, that projected paired kernels
are nearly $S^*$-invariant subspaces of $H^2_+$, and so are their closures, by the following 
result.

\begin{prop}\label{prop:n4.4}
If $\GS \subset H^2_+$ is nearly $S^*$-invariant,
then its closure is also nearly $S^*$-invariant.
\end{prop}
\beginpf
Clearly we may assume without loss of generality that $\GS \ne \{0\}$.
In that case, since $\GS$ is nearly $S^*$-invariant, there must exist $h \in \GS$ with $h(0)=1$.
Now if $f \in \bar\GS$ with $f(0)=0$, then
there is a sequence $(f_n)$ in $\GS$ with $f_n \to f$ in norm and
hence $f_n(0) \to f(0)=0$.
Thus for each $n$ we have that $f_n-f_n(0)h$ is a function in $\GS$ vanishing at $0$,
so $f_n-f_n(0)h=zg_n$ for some
$g_n \in \GS$, and $\lim f_n=f = \lim zg_n$.
Hence $(g_n)$ converges to a function $g \in \bar\GS$ such that $f=zg$.
\endpf

The closed nearly $S^*$-invariant subspaces of $H^2_+$ admit
a representation as a product of the form $uK_\theta$, where $u \in H^2_+$ and
$K_\theta$ is a model space, which will be considered in the
next section. In the case of
the closure of $\k^+_{a,b}$, $u$ must be outer by Corollary \ref{ncor:5.3}.

The notion of near $\eta$-invariance can naturally be extended to $\k^-_{a,b}$, replacing
$H^2_+$ by $H^2_-$. The following
proposition
shows that it is enough to consider the problem of near invariance for $\kerab^+$.
    
     \begin{prop}\label{nprop:5.1}
  $\kerab^-$ is nearly $\bar\eta$-invariant in $H^2_-$ if and only if
  $\k^+_{\bar b,\bar a}$ is nearly $\eta$-invariant in $H^2_+$.
  \end{prop}
  
  \beginpf
  This is a consequence of Proposition \ref{nprop:2.4}.
  Suppose that $\kerab^-$ is nearly $\bar\eta$-invariant in $H^2_-$, i.e.,
  \[
  \phi_- \in \k^-_{a,b}, \quad \bar\eta \phi_- \in H^2_- \implies \bar\eta \phi_- \in \k^-_{a,b},
  \]
  and let
  \[
  \phi_+ \in \ker^+_{\bar b,\bar a}, \qquad \eta \phi_+ \in H^2_+.
  \]
  Then $\bar z \ol{\phi_+} \in \kerab^-$, $\bar\eta \bar z \ol{\phi_+} \in H^2_-$,
  so $\bar\eta \bar z\ol{\phi_+} \in \kerab^-$ by near invariance which,
by Proposition \ref{nprop:2.4}, implies that $\eta \phi_+ \in \k^+_{\bar b,\bar a}$. So $\k^+_{\bar b,\bar a}$ is nearly $\eta$-invariant in $H^2_+$.
  The converse is proved analogously.  
\endpf

Since $\kerab^+$ is nearly $\bar z$-invariant in $H^2_+$ (equivalently, 
nearly $S^*$-invariant), it follows that $\kerab^-$ is nearly $z$-invariant in $H^2_-$. Therefore,
if $\kerab^\pm \ne \{0\}$, there exists $\phi_+ \in \kerab^+$ with $\phi_+(0) \ne 0$
and there exists $\phi_- \in \kerab^-$ such that $(z\phi_-)(\infty) \ne 0$
(i.e., $ \ol{\psi_+(0)} \ne 0$,
where $\psi_+= \bar z \ol{\phi_-}$).
Analogously, since $\kerab^+$ is nearly $\frac{1}{z-z_0}$-invariant for any $z_0 \in \DD$, there exists, for any $z_0 \in \DD$,  $\phi_+ \in \kerab^+$ with $\phi_+(z_0) \ne 0$
and $\psi_- \in \kerab^-$ such that $\psi_- (1/z_0) \ne 0$.

\section{Minimal Toeplitz kernels and representations of projected paired kernels}
\label{sec:n5}

If $a/b \in L^\infty$, then $\k^+_{a,b}$ is a Toeplitz kernel; but, in general, $\k^+_{a,b}$ may not even be a closed subspace
of $H^2_+$, as in the case studied in Section \ref{sec:3}. There it was also shown that,
although $\k^+_{\bar\theta,z+1}$ is not a Toeplitz
kernel, one can nevertheless  determine a minimal Toeplitz kernel containing it. It is thus natural
to ask if such a property holds for every projected paired kernel. The answer is in the
affirmative, as one of the consequences of the following theorem.

Recall that a {\em maximal function\/} $\phi_m$ for a Toeplitz kernel $\ker T_g$
is one such that $\ker T_a$ is the minimal kernel to which $\phi_m$ belongs (see \cite{CP14}).
Every Toeplitz kernel possesses a maximal function.

\begin{thm}\label{thm:5.1bis}
Let  $a \in L^\infty \setminus \{0\}$ and $b \in H^2_+$ such that $\ker T_a \ne \{0\}$
and $b \ker T_a \subset H^2_+$.
Then there exists a minimal Toeplitz kernel containing $b \ker T_a$, which is 
$\ker T_{a\bar b/b_o}$, where
$b_o$ is the outer factor of $b$.
Moreover, if $\phi_m$ is a maximal function for $\ker T_a$, then $b \phi_m$
is a maximal function for $\ker  T_{a\bar b/b_o}$.
\end{thm}

\beginpf
If $\ker T_a \ne \{0\}$, let $\phi_m$ be a maximal function of
$\ker T_a$, and write
$\phi_m=I_+O_+$ with
$I_+$ inner and $O_+$ outer in $H^2_+$.
Then 
\beq \label{neq:n5.1}
\ker T_a= \ker T_{\bar z\overline{I_+}\overline{O_+}/O_+}
\eeq
\cite[Thm. 5.1]{CP14}. 
On the other hand, there exists a minimal
kernel for $b \phi_m \in b \ker T_a \subset H^2_+$ which,
if $b=b_i b_o$ is an inner--outer factorization of $b$ with $b_i$ inner and $b_o$ outer, is given by
\beq\label{neq:4.1}
\Kmin(b\phi_m) = \ker T_{
\dfrac
{\bar z\overline{I_+}\overline{b_i} \overline{b_o} \overline{O_+}}
{b_o O_+ }
}
= \ker T_{
\dfrac{\bar z\overline{b \phi_m}}{b_o O_+ }}.
\eeq
Let us show that 
\beq\label{neq:4.2}
\Kmin(b \phi_m) \supset b \ker T_a.
\eeq
Let $\psi_+$ be any non-zero element of $\ker T_a$;
then (see \eqref{neq:n5.1}) $\dfrac{\bar z\overline{I_+}\overline{O_+}}{O_+} \psi_+ \in H^2_-$ and
\beq
\underbrace{\frac
{\bar z\overline{I_+}\overline{O_+}\overline{b_i} \overline{b_o} }
{O_+ b_o}}_{\in L^\infty} 
\underbrace{(b \psi_+)}_{\in H^2_+} = 
\underbrace{\overline{b_o}}_{\in \ol{H^2_+}}
\underbrace{\frac{\bar z\overline{I_+}\overline{O_+}}{O_+} \psi_+}_{\in H^2_-=\bar z \ol{H^2_+}} \in \bar z \ol{H^1} \cap L^2 \subset \bar z H^2_+=H^2_-.
\eeq
Thus $b\psi_+ \in \Kmin(b \phi_m)$ as defined in \eqref{neq:4.1}, and \eqref{neq:4.2} holds.
On the other hand, since $b\phi_m \in b \ker T_a$, any Toeplitz kernel
containing $b \ker T_a$ must also contain $\Kmin(b \phi_m)$,
so the latter is the minimal kernel containing $b \ker T_a$.
Finally, from $\ker T_a = \ker T_{\bar z\overline{I_+}\overline{O_+}/O_+}$
we conclude that 
\[
a= \dfrac{\bar z\overline{I_+}\overline{O_+}}{O_+} h_- = \dfrac{\bar z\overline{\phi_m}}{O_+} h_-
\]
for some $h_- \in \G \overline{H^\infty}$ (by \cite[Cor. 7.8]{CP18b}. Therefore,
\[
a \frac{\bar b}{b_o} = \frac{\bar z\overline{b \phi_m}}{b_o O_+} h_-
\]
and it follows from \eqref{neq:4.1} that $ \Kmin (b \phi_+)=\ker T_{a \bar b/b_o}$.

\endpf

\begin{cor}\label{cor:n5.2}
Every   subspace of $H^2_+$ of the form $uK_\theta$,
where $u \in H^2_+$ and $\theta$ is an inner function, is contained
in a minimal Toeplitz kernel.
\end{cor}

As an example, take the example of Section \ref{sec:3}, namely $(z+1)K_\theta$. By Theorem 
\ref{thm:5.1bis} we have that the minimal kernel containing 
$(z+1)K_\theta$ is 
\[
 \ker T_{\bar\theta \frac{\ol{z+1}}{z+1}}= \ker T_{\bar\theta \bar z}=K_{\theta z},
 \]
 as shown before.\\

A well-known theorem by Hitt \cite{hitt} describes the closed nearly $S^*$-invariant subspaces of $H^2_+$ as having
the form $M=uK$, where $u \in H^2_+$ has unit norm, $u(0)>0$,
$u$ is orthogonal to all elements of $M$ vanishing at the origin, $K$ is
an $S^*$-invariant subspace, and the operator of multiplication by $u$
is isometric from $K$ into $M$. Naturally, one can have
$K=\{0\}$ or $K=H^2_+$, but the most interesting cases are those in which $K$ is a model space $K_\theta=\ker T_{\bar\theta}$.

\begin{cor}
For every nondegenerate $(a,b)$ there exists a minimal kernel containing $\k^+_{a,b}$,
which coincides with $\k^+_{a,b}$ if $a/b \in L^\infty$.
\end{cor}
\beginpf
By Proposition \ref{prop:n4.4}, the closure of $\kerab^+$ is nearly $S^*$-invariant, so we have that
the closure of $\kerab^+$ is $uK$, where $u \in H^2_+$ is outer, by Corollary \ref{ncor:5.3},
and $K=H^2_+$ or $K=K_\theta$ with $\theta$ inner. In the 
latter case, the
result follows from 
Proposition \ref{prop:n4.1}, Proposition \ref{prop:n4.4}, Hitt's theorem, and 
Corollary \ref{cor:n5.2}.
If $K=H^2_+$, then, since $uH^2_+ \subseteq H^2_+$ is closed and $u$ is outer, we must have 
$u \in \G H^\infty$ and $uK=H^2_+$.
\endpf

Hayashi showed in \cite{hayashi} that the kernel of every Toeplitz operator $T_g$  can be
written as $uK_\theta$ where $u$ is outer, $u^2$ is rigid (an exposed point of the unit ball of $H^1$), $\theta$ is inner with
$\theta(0)=0$, and $u$ multiplies $K_\theta$ isometrically onto the Toepliz kernel.
It may happen that $u \in \G H^\infty$ and, in that case, we can write
\beq
uK_\theta = u \ker T_{\bar\theta} = \ker T_{\bar \theta u^{-1}} = \k^+_{\bar\theta,u}.
\eeq

Other representations of a similar form can be found for Toeplitz
kernels. For instance, if
$g \in L^\infty$ admits a Wiener--Hopf factorization of the form $g=g_- \bar\theta g_+$ with
$g_- \in \G \ol{H^\infty}$, $g_+ \in \G H^\infty$ and $\theta$ inner, then
\beq
\ker T_g = g_+^{-1} \ker T_{\bar\theta} = g^{-1}_+ K_\theta = \k^+_{\bar\theta,g^{-1}_+}.
\eeq
Another example where
\beq\label{eq:n5.7}
\k^+_{a,b}=b \ker T_a
\eeq
is the case studied in Section \ref{sec:3}, which is not a Toeplitz kernel. It is thus natural to ask for conditions under which
\eqref{eq:n5.7} holds and, in general, what is the relation between $\k^+_{a,b}$ and $b \ker T_a$ for $a,b \in L^\infty$
nondegenerate.
We have the following.

\begin{prop} \label{prop:n5.4}
For $a,b \in L^\infty$ we have that
$ b\ker T_a \cap H^2_+ \subseteq \kerab^+$ and
\beq\label{eq:n5.8}
\kerab^+=b\ker T_a \cap H^2_+ \,\,\,{\rm if \,and\,  only \, if}\,\,\,
 \kerab^+ \subseteq bH^2_+.
\eeq
\end{prop}

\beginpf
We have $a \ker T_a \subseteq H^2_-$, so, for any $\phi_+ \in \ker T_a$ such that
$b\phi_+ \in H^2_+$,
\[
a \underbrace{(b\phi_+)}_{\in H^2_+} + b\underbrace{(-a\phi_+)}_{\in H^2_-} = 0,
\]
and it follows that $b\phi_+ \in \kerab^+$. So $ b\ker T_a \cap H^2_+ \subseteq \kerab^+$.

 Now, it is clear that $b\ker T_a \cap H^2_+=\kerab^+$ implies that
$\kerab^+ \subseteq bH^2_+$.
Conversely, suppose that $\kerab^+ \subseteq bH^2_+$.
Then, for any $\phi_+ \in \kerab^+$ we have $\phi_+ =b\psi_+$
with $\psi_+ \in H^2_+$, and therefore, for some $\phi_-\in H^2_-$,
\begin{eqnarray*}
a\phi_+ + b \phi_- = 0 & \iff & ab \psi_+ + b \phi_- = 0 \\
&\iff& b(a\psi_+ + \phi_-) = 0.
\end{eqnarray*}
We have the standing assumption  that $b \ne 0$ a.e.\ on $\TT$, so $a\psi_+ = -\phi_-$.
We conclude that $\psi_+ \in \ker T_a$, so $\phi_+ \in b\ker T_a \cap H^2_+$, and it
follows that
$\kerab^+ \subseteq b \ker T_a \cap H^2_+.
$
\endpf

An immediate consequence of Proposition \ref{prop:n5.4} is the following.

\begin{cor}
If $b \in H^\infty$ then $b \ker T_a \subseteq \k^+_{a,b}$.
\end{cor}

It is clear from the invariance results of Section \ref{sec:n4} that we can have $\k^+_{a,b} \subseteq bH^2_+$, with $b \in H^\infty$, only if $b$ is outer.
Moreover, we have the following.

\begin{cor}\label{cor:n5.6}
If $b \in H^\infty$ is outer and either (i) $a \in \G L^\infty$ or
(ii) $b \in \G H^\infty$, then $\k^+_{a,b}=b \ker T_a$.
\end{cor}

\beginpf
If $b \in \G H^\infty$ then $bH^2_+=H^2_+$ and the equality follows from \eqref{eq:n5.8}.
If $a \in \G H^\infty$ then
\[
a\phi_+ + b \phi_- = 0 \iff 
\underbrace{\frac{\phi_+}{b}}_{\in \GN^+}
= \underbrace{-\frac{\phi_-}{a}}_{\in L^2}
\in \GN^+ \cap L^2 = H^2_+,
\]
so $\phi_+ \in b H^2_+$ and again the equality in Corollary \ref{cor:n5.6}
follows from \eqref{eq:n5.8}.
\endpf

\section{Inclusion relations} \label{sec:6}

As in the case of Toeplitz kernels \cite{CP14}, the near invariance properties of projected paired kernels imply certain lower bounds for the dimension of a paired kernel
containing a given function. For instance, if
$\phi \in \kerab$ and $\phi_+(0)=0$, then $\kerab$ must also
contain the function $\psi \in L^2$ with $\psi_+=\bar z \phi_+$ (note that this
defines $\psi$ by Proposition \ref{nprop:4.2}). 
As another example, using a similar reasoning, if there exists $\phi\in \kerab$ such that
$\phi_+ \in \theta H^2_+$ or $\phi_- \in \bar\theta H^2_-$ where $\theta$
is inner but not a finite Blaschke product, then $\kerab$ is infinite-dimensional.

On the other hand, it is easy to see that, if $\theta_1$ and $\theta_2$ are
inner functions, then
\beq\label{neq:6.1}
\k^+_{a\theta_1,b \ol{\theta_2}} \subseteq \kerab^+,
\eeq
We may then ask if the inclusion is strict and,
in that case, how much ``smaller'' 
$\k^+_{a\theta_1,b \ol{\theta_2}}$ is with respect to $\kerab^+$ and,
in particular, when it is $\{0\}$.
More generally, one may ask what are the relations between two paired kernels, or two projected paired kernels, 
whose symbol pairs are related by multiplication operators.

Note  that, since  a nontrivial paired kernel cannot be contained in a different one, and indeed their intersection is $\{0\}$ (see \cite{CGP23}), obtaining, for example,
 a paired operator analogue of the property $\ker T_{\theta g} \subsetneq \ker T_g$ (valid for
$g \in L^\infty$ and $\theta$ inner, nonconstant) is possible only
by saying that $\ker_{\theta a,b}$ is isomorphic to a proper subspace of $\ker _{a,b}$.
Alternatively, taking Proposition  \ref{nprop:4.2}   into account, we can say it in an equivalent and simpler way as $\k^+_{\theta a,b} \subsetneq \k^+_{a,b}$ (cf. Proposition
 \ref{nprop:6.1}).

In the following propositions we present several inclusion relations between projected paired kernels which generalise similar properties valid for Toeplitz kernels,
using in particular the near invariance properties of
Section \ref{sec:n4} to establish strict inclusions.

We recall that, for any $\eta \in L^\infty$,
$\kerab=\k_{a\eta,b\eta}$, so $\k^+_{a,b}=\k^+_{a\eta,b\eta}$.

\begin{prop}\label{nprop:6.1}
(i) If $h_- \in \ol{H^\infty}$ then $ \k^+_{a,bh_-}
\subseteq \kerab^+ \subseteq \k^+_{ah_-,b}$.\\
(ii) If $\ol{h_-}$ is outer, then
\begin{eqnarray*}
(a) \quad \k^+_{a,bh_-}=\kerab^+ && \hbox{if } \frac{a}{bh_-} \in L^\infty, \\
(b) \quad \k^+_{ah_-,b}=\kerab^+ && \hbox{if }\frac{a}{b} \in L^\infty.
\end{eqnarray*}
(iii) If $\ol{h_-}$ has a non-constant inner factor, then
\[
\k^+_{a,bh_- } \subsetneq \k^+_{a,b} \subsetneq \k^+_{ah_-,b}.
\]
\end{prop}
\beginpf
(i) $a\phi_++bh_- \phi_- = 0 \implies a\phi_++b (h_-\phi_-)=0$ and\\
$a\phi_++b\phi_-=0 \implies (ah_-)\phi_+ + b (h_- \phi_-) = 0$.\\

(ii) Let now $\ol{h_-}$ be outer. We start by proving the second equality
(ii)(b). We have
\[
ah_- \phi_+ + b \phi_- = 0 
\iff \dfrac{\phi_-}{h_-} = - \dfrac{a}{b} \phi_+
\iff \dfrac{ \ol{\phi_-}}{\ol{h_-}} = - \ol{\dfrac{a}{b} \phi_+}.
\]
Since the left-hand side of the last equation is in the Smirnov class $\GN_+$
when $\ol{h_-}$ is outer and the right-hand side is in $L^2$ if
$a/b \in L^\infty$, we have under these assumptions that $\phi_-/h_- \in H^2_-$,
so $\phi_+ \in \kerab^+$ since $a\phi_+ + b \dfrac{\phi_-}{h_-}= 0$, and thus 
$\k^+_{ah_-,b} \subseteq \kerab^+ $. The equality follows from (i). \\

Next, we   have that
\[
a\phi_++b\phi_-=0 \implies a\phi_+ + bh_- \dfrac{\phi_-}{h_-} =0
\implies \dfrac{\phi_-}{h_-} = - \dfrac{a}{b h_-} \phi_+
,
\]
and we conclude analogously that $\dfrac{\phi_-}{h_-} \in H^2_-$ so
$\kerab^+ \subseteq \k^+_{a,bh_-}$. Again the equality (ii)(a) follows from (i).\\

(iii) Suppose that $\ol{h_-}$ has a non-constant inner factor.
If $\kerab^+ \subseteq \k^+_{a,bh_-}$ then, for any $\phi \in L^2$ such that
$a\phi_+ + b\phi_-=0$ we must also have $\psi_-  \in H^2_-$ such that
$a\phi_+ + bh_- \psi_-=0$.
Thus $\phi_-=h_- \psi_-$. This implies that $\kerab^- \subseteq h_- H^2_-$, which is impossible by Corollary \ref{ncor:5.3} because $\ol{h_-}$ has a nonconstant inner factor. So
$\k^+_{a,bh_-} \subsetneq \kerab^+$.

A similar argument shows that $\k^+_{a,b} \subsetneq \k^+_{ah_-,b}$.
\endpf


\begin{prop}
Suppose that $a/b \in L^\infty$ and $z_0 \in \TT$. Then\\
$\k^+_{a(z-z_0),b}=\k^+_{az,b}$.
\end{prop}
\beginpf
$\k^+_{a(z-z_0),b} = \k^+_{az (z-z_0)/z,b}= \k^+_{az,b}$
by Proposition \ref{nprop:6.1}(ii), since $(z-z_0)/z \in \ol{H^\infty}$ and
its conjugate is $1-\ol{z_0}z$, an outer function in $H^\infty$.
\endpf

As an immediate corollary we have:

\begin{cor}
Let $g \in L^\infty$ and $z_0 \in \TT$. Then
$\ker T_{g(z-z_0)^n}=\ker T_{g z^n}$.
\end{cor}

\begin{rem}{\rm
This corollary allows us to generalize several results from \cite[Sec. 6]{CMP16},
in particular Theorems 6.2 and 6.7 in \cite{CMP16}, which deal
with the relations between $\ker T_g$ and $\ker T_{\theta g}$, for a
finite Blaschke product $\theta$, to the case where the symbol has
zeros of integer order on $\TT$.}
\end{rem}

\begin{prop}\label{nprop:6.5}
(i) If $h_+ \in H^\infty$ then
\[
h_+ \k^+_{ah_+,b} \subseteq \kerab^+ \quad \hbox{and} \quad
h_+ \kerab^+ \subseteq \k^+_{a,bh_+}.
\]
(ii) If $h_+ \in \G H^\infty$, then
\[
h_+ \k^+_{ah_+,b} = \kerab^+ = h_+^{-1} \k^+_{a,bh_+}.
\]
(iii) If $h_+$ has a non-constant inner factor, then
\[
h_+ \k^+_{ah_+,b} \subsetneq \kerab^+ 
\quad \hbox{and} \quad
h_+ \kerab^+ \subsetneq \k^+_{a,bh_+}.
\]
\end{prop}
\beginpf
(i) We can write $(ah_+)\phi_++b\phi_- = 0$ as $a(h_+\phi_+) +b\phi_- = 0$,
from which $h_+ \k^+_{ah_+,b} \subseteq \kerab^+$; then if $\phi_+ \in \kerab^+$ we have
$a\phi_++b\phi_-=0$ for some $\phi_- \in H^2_-$, so,
since $a(h_+ \phi_+)+(b h_+)\phi_-=0$, we have $h_+\phi_+ \in \k^+_{a, bh_+}$.\\
(ii) If $h_+ \in \G H^\infty$, then from $a\phi_++b\phi_-=0$ we obtain
$ah_+(\phi_+ h_+^{-1})+b \phi_-=0$, so
$\kerab^+ \subseteq h_+ \k^+_{ah_+,b}$ and we have equality from (i).
The second equality follows from $h_+^{-1}\k^+_{a h_+^{-1},b} = h_+^{-1}\k^+_{a,bh_+}$.\\
(iii) If $h_+$ has a non-constant inner factor, we cannot have 
$\kerab^+ \subseteq h_+ H^2_+$ nor $\k^+_{a,bh_+}\subseteq h_+ H^2_+$ by the near-invariance result of Corollary \ref{ncor:5.3}.
\endpf

As a consequence of Proposition \ref{nprop:6.1} (ii) and Proposition \ref{nprop:6.5} (ii) we have
the following.

\begin{prop}\label{nprop:6.6}
Let $B=B_- g B_+ \in L^\infty$, where
$B_-^{\pm 1} \in \ol{H^\infty}$
and $B_+^{\pm 1} \in H^\infty$. Then
\[
\k^+_{aB,b} = B_+^{-1} \k^+_{ag,b} \quad \hbox{and} \quad \k^+_{a,Bb}=B_+ \k^+_{a,gb}.
\]
\end{prop}

\beginpf
\[
\phi_+ \in \k^+_{aB,b} \iff \exists \phi_-:  aB\phi_++b\phi_-=0,
\]
that is,
$aB_-gB_+\phi_+ + b\phi_-=0$.\\
We write this as $\exists \phi_- :   ag(B_+\phi_+)+b(\phi_-B^{-1}_-)=0$.
Equivalently,
$\exists \psi_- :   ag(B_+\phi_+)+b\psi_-=0$.\\
Finally, this holds if and only if $B_+\phi_+ \in \k^+_{ag,b}$, i.e., 
$\phi_+ \in B^{-1}_+ \k^+_{ag,b}$.

The other identity is proved similarly.
\endpf


\begin{cor}\label{cor6.6a}

If $h_- \in \G \ol{H^\infty}$ then
$\kerab^+ = \k^+_{ah_-,b} = \k^+_{a,h^{-1}_-b} = \k^+_{a,h_-,b}$.

\end{cor}

\beginpf
We have $\kerab^+ \subseteq \k^+_{ah_-  ,b}$; conversely,
\begin{eqnarray*}
ah_- \phi_+ + b \phi_- = 0 & \implies & a\phi_+ + b (h^{-1}_- \phi_-) = 0 \\
& \implies & \phi_+ \in \kerab^+,
\end{eqnarray*}
so $\kerab^+ = \k^+_{ah_-,b}= \k^+_{a,bh^{-1}_- }$ and,
replacing $h_-$ by $h^{-1}_-$ we conclude that
$\kerab^+ = \k^+_{a,bh_-}$.

\endpf


Since the inclusion in \eqref{neq:6.1} is strict when $\theta_1$ or $\theta_2$ is
not constant or, equivalently, 
$\k^+_{a\theta,b} \subsetneq \k^+_{a,b}$ when $\theta$ is a non-constant inner function, one
may ask whether the dimensions of 
those two spaces can be compared, as in the case of Toeplitz kernels.
 The following
theorem generalises analogous results obtained in \cite[Sec. 2]{BCD10}
and \cite[Sec. 6]{CMP16} for Toeplitz kernels, and can be proved
in a similar way.

\begin{thm}\label{nthm:6.7}
Let $\theta$ be a non-constant finite Blaschke product. Then\\
(i) $\k^+_{a\theta,b}$ is finite-dimensional if and only if $\kerab^+$ is finite-dimensional.
Similarly for $\k^+_{a,\theta b}$.\\
(ii) If $\dim \kerab^+ = d < \infty$ and $\dim K_\theta \ge d$ then $\k^+_{a\theta,b}=\{0\}$.\\
(iii) If $\dim \kerab^+ = d < \infty$ and $\dim K_\theta = k < d$, then
$\dim\k^+_{a\theta,b}=d-k$.\\

\end{thm}

Note that by Proposition \ref{nprop:4.2} and Corollary \ref{cor:2.2} the same result holds if we replace
$\kerab^+$ by $\kerab$ etc.
When $\kerab^+$ is not finite-dimensional, one
can still compare it with $\k^+_{a\theta,b}$ by means of  
the following decomposition.

\begin{thm}
Let $\theta$ be a non-constant finite Blaschke product.
If $\theta$ has $k$ zeros in $\DD$, counting multiplicities,
in which case we can write $\theta=B_- z^k B_+$ with $B_-^{\pm 1} \in \ol{H_\infty}$ and
$B_+^{\pm 1} \in H^\infty$, then
there exist $\psi_{j+} \in H^2_+$ ($j=0,\ldots,k-1$) with
$\psi_{j+}(0)=1$, for each $j$ such that the following direct sum decomposition holds:
\begin{eqnarray*}
\kerab^+ &=& z^k B_+ \k^+_{a\theta,b}
+ \spam \{\psi_{0+},z\psi_{1+},\ldots, z^{k-1} \psi_{(k-1)+} \}\\
&=& z^k \k^+_{az^k,b} + \spam \{\psi_{0+},z\psi_{1+},\ldots, z^{k-1} \psi_{(k-1)+} \}.
\end{eqnarray*}
\end{thm}
\beginpf
By Proposition \ref{nprop:6.6} it is enough to consider
$\k^+_{az^k,b}$.
Since $\kerab^+$ is nearly $S^*$-invariant,
there exists $\psi_{0+} \in \kerab^+$ with
$\psi_{0+}(0)=1$;
let $\psi_{0-}$ be given by  $a\psi_{0+} + b \psi_{0-} =0$. Then for any $\phi_+ \in \kerab^+$,
\begin{eqnarray*}
a\phi_+ + b \phi_- = 0 & \iff & az \frac{\phi_+-\phi_+(0)\psi_{0+}}{z} + a \phi_+(0) \psi_{0+}+b\phi_-=0\\
& \iff &  az \frac{\phi_+-\phi_+(0)\psi_{0+}}{z} + b(\phi_- - \phi_+(0)\psi_{0-})=0,
\end{eqnarray*}
so 
\[
\tilde \phi_+ := \frac{\phi_+-\phi_+(0)\psi_{0+}}{z} \in \k^+_{az,b}
\]
 and
$\phi_+ = z\tilde\phi_+ + \phi_+(0) \psi_{0+}$.

 Therefore $\kerab^+ = z\k^+_{az,b} \oplus \spam\{\psi_0\}$.
 
 Proceeding analogously with $\k^+_{az,b}, \k^+_{az^2,b},\ldots,\k^+_{az^{k-1},b}$,
 we get 
 \[
 \kerab^+ = z^k \k^+_{az^k,b}
 \oplus \spam \{\psi_{0+},z\psi_{1+},\ldots,z^{k-1}\psi_{(k-1)+}\},
 \]
 where $\psi_{j +} \in H^2$ and $\psi_{j+}(0)=1$ for all $j=0,1,\ldots,k-1$. 
\endpf 

We now obtain a description of some related
projected paired kernels that will be used in the next section
to study truncated Toeplitz operators: $\k^+_{p_1,p_2}$, $\k^+_{\alpha p_1,p_2}$ and
$\k^+_{p_1,\alpha p_2}=\k^+_{\bar\alpha p_1,p_2}$, where $p_1$ and $p_2$
are polynomials of degrees $n_1$ and $n_2$ respectively
 without common zeros, and $\alpha$ is an inner function.
 Note that $\k^+_{p_1,p_2}$ can be related to a new class of Toeplitz-like
 operators introduced in 
 \cite{GHJR1,GHJR2,GHJR3}: see \cite{CGP23}.

\begin{prop}\label{nprop:6.8}
Let $p_1$ and $p_2$ satisfy the assumptions above,
and let moreover $p_i=p_{i\DD}p_{i\TT}p_{i\EE}$, for $i=1,2$, where
the zeros of $p_{i\DD}$, $p_{i\TT}$ and $p_{i\EE}$
are in $\DD,\TT$ and $\EE=\DD^e$, respectively, and we denote
by $n_{i\DD},n_{i\TT}$ and $n_{i\EE}$ the corresponding degrees. Then\\
(i) $\k^+_{p_1,p_2}=\{0\}$ if $n_2 \le m:= n_{2\TT}+n_{1\TT}+n_{2\EE}+n_{1 \DD}$;\\
(ii) $\k^+_{p_1,p_2} = p_{1\TT}p_{2\TT} p_{1\DD}p_{2\EE} \GP_{n_2-m-1}$,
where $\GP_\ell$, for $l\in \NN \cup \{0\}$, denotes the space of all
polynomials of degree less than or equal to $\ell$, and
\[
\dim \k^+_{p_1,p_2} = n_2-m= n_{2\DD}-n_{1\DD}-n_{1 \TT}
\]
if $n_2>m$.
\end{prop}
\beginpf
\[
p_1\phi_+ + p_2 \phi_- = 0
\iff
p_1 \phi_+ = -p_2 \phi_-
= -\underbrace{p_{2\TT}p_{2\EE}p_{2\DD}}_{\text{degree } n_2} \phi_-
= q_{n_2-1},
\]
where $q_{n_2-1}$ is a polynomial of degree less than or equal to $n_2-1$.
Since
$\dfrac{q_{n_2-1}}{p_2} \in H^2_-$ and $\dfrac{q_{n_2-1}}{p_1} \in H^2_+$,
$q_{n_2-1}$ must be of the form
$q_{n_2-1}=p_{2\TT}p_{2\EE}p_{1\TT}p_{1\DD}q$, where $q$ is a polynomial.
Thus if $n_2-1 < m$ then $q_{n_2-1}=0$ and (i) holds.
Also, if $n_2>m$ then (ii) holds.
\endpf

\begin{prop}\label{nprop:6.9}
With the same assumptions and notation as in Proposition \ref{nprop:6.8}, let $\alpha$ be a nonconstant inner function. In the case that
$\alpha$ is a finite Blaschke product,
let $\deg\alpha$ denote the number of zeros of $\alpha$, counting multiplicity. 
Let $d=n_{2\DD}-n_{1\DD}-n_{1\TT} > 0$.\\
(i) $\k^+_{\alpha p_1,p_2}=\{0\}$ if $\alpha$ is not a finite Blaschke product or
$\deg\alpha \ge d > 0$.\\
(ii) $\dim\k^+_{\alpha p_1,p_2} = d -\deg \alpha$ if $\alpha$ is a finite Blaschke product with
$\deg\alpha<d$.
In this case, factorizing $\alpha=B_- z^{\deg\alpha} B_+$ with
$B_-^{\pm 1} \in \ol{H^\infty}$ and $B_+^{\pm 1} \in H^\infty$, we have that
$\k^+_{\alpha p_1,p_2} = B^{-1}_+ \k^+_{z^{\deg\alpha}p_1,p_2}$,
where the projected paired kernel on the right-hand side is described in Proposition
\ref{nprop:6.8}
\end{prop}
\beginpf
(i) If $\alpha$ is not a finite Blaschke product or $\deg\alpha \ge d$, then (i) holds by Theorem \ref{nthm:6.7}.\\
(ii) If $\deg\alpha < d$ then 
\[
\k^+_{\alpha p_1,p_2}=\k^+_{B_- z^{\deg\alpha}B_+p_1,p_2}=B^{-1}_+ \k^+_{z^{\deg\alpha}p_1,p_2}
\]
by Proposition \ref{nprop:6.6}.
\endpf
Analogously, we have:

\begin{prop}\label{prop:6.12bis}
With the same assumptions and notation as in Proposition \ref{nprop:6.9},
and assuming that $\dim\k^+_{p_1,p_2}=d=n_{2\DD}-n_{1\DD}-n_{1\TT}>0$, we have
\[
\k^+_{\bar\alpha p_1,p_2} = \ol{B_-^{-1}} \k^+_{\bar z^{\deg\alpha} p_1,p_2}
= \ol{B_-^{-1}} \k^+_{p_1,z^{\deg\alpha}p_2},
\]
where $\alpha=B_- z^{\deg\alpha}B_+$ as in Proposition \ref{nprop:6.9} and
$\k^+_{p_1,z^{\deg\alpha}p_2}$ is described in Proposition \ref{nprop:6.8}.
\end{prop}

\beginpf
The first equality follows from Proposition \ref{nprop:6.6} and the second from
the definition of $\k^+$.
\endpf

\section{Kernels of finite-rank asymmetric truncated Toeplitz operators}
\label{sec:7}

Paired operators can also be defined in the matricial setting. They appear in the
literature \cite{MP,CG,CKP}
as operators on $(L^2)^n$ with $n \times n$ matricial coefficients $A,B \in \GL((L^2)^n)$.

In most cases $A$ and $B$ are multiplication operators and in that case the paired
operator takes the form
\beq
T_{A,B}=AP^+ + BP^- \qquad \hbox{with} \quad A,B \in (L^\infty)^{n \times n}.
\eeq
One can, however, consider other possible generalizations of scalar paired operators, for instance
\beq
T_{A,B}: (L^2)^n \to L^2, \qquad T_{A,B}=AP^++BP^- \quad \hbox{with} \quad A,B \in \GL((L^2)^n,L^2)
\eeq
or, considering in particular multiplication operators,
\beq
T_{A,B}=AP^++BP^-: (L^2)^n \to L^2 \qquad \hbox{with} \quad A,B \in (L^\infty)^{1 \times n},
\eeq
\[
A=[a_i]^T_{i=1,\ldots,n}, \qquad B=[b_i]^T_{i=1,\ldots,n}.
\]
The kernels of operators of this form are defined by the Riemann--Hilbert problem
\beq
\sum_{i=1}^n a_i \phi_{i+} = - \sum_{i=1}^n b_i \phi_{i-}, \qquad \phi_{i\pm} \in H^2_{\pm}.
\eeq
The latter have appeared in recent works: see, for instance,
\cite{ACCR} and \cite{CP20}.\\

We apply here the results of the latter, together with those of previous sections,
to study the kernels of finite-rank asymmetric truncated Toeplitz operators (ATTO for short) of the form
\beq\label{eq:7.5}
A_\phi^{\theta,\alpha}: K_\theta \to K_\alpha, \qquad A_\phi^{\theta,\alpha}=P_\alpha \phi {P_\theta}_{|K_\theta},
\eeq
where $\theta,\alpha$ are inner functions, $P_\theta$ and $P_\alpha$
denote the orthogonal projections from $L^2$ onto $K_\theta$ and $K_\alpha$ respectively,
and $\phi \in L^\infty$ is called the {\em symbol\/} of the ATTO.

Indeed, if we define $P_1$ to be the projection $P_1(x,y)=x$,
then we have
\beq
\ker A_\phi^{\theta,\alpha}= P_1 \ker T_G,
\eeq
where $G$ has the matrix symbol
\beq \label{eq:7.6bis}
G=\left[
\begin{matrix}
\bar\theta & 0 \\ \phi & \alpha \end{matrix}
\right].
\eeq
Now it was shown in \cite[Thm. 3.1 and Cor. 3.4]{CP20}  that the kernels of block
Toeplitz operators with symbols of the form \eqref{eq:7.6bis}
are of so-called {\em scalar type}, i.e., can be described as
scalar multiples of a fixed vector function, and can be
obtained from a pair of functions  $f,g \in \GF^2$, 
where $\GF$ consists of all complex-valued functions defined a.e.\ on $\TT$,
provided that $f$ and $g$ satisfy the relation
$Gf=g$.

To present the results from \cite{CP20} that will be used here, it is useful
to settle some notation, as follows:

(i) if $f=(f_1,f_2) \in \GF^2$, then we say that $f$ is
left-invertible in $\GF$ if and only if there exists $\tilde f \in \GF^2$ such that
$\tilde f^T f=1$, and, in that case, $\tilde f^T$ is called the left
inverse of $f$ in $\GF$;

(ii) we define $\GH_X = \{h \in \GH: Xh \in \GH \}$
for any subspace $\GH \subseteq (L^2)^n$ and any $X \in \GF^{n \times n}$;

(iii) for any matricial coefficients $A,B$ with elements in $\GF$ let
us define $T_{A,B} f = AP^+f + BP^- f$ for
$f \in (L^2)^n$ and 
write ${\frak s}_{A,B}=A(H^2_+)^n \cap B(H^2_-)^n$.\\

With this notation we can formulate Corollary 3.4 in
\cite{CP20} as follows.

\begin{thm}\label{thm:7.1}
Let $f=(f_1,f_2)$ and $g=(g_1,g_2)$ belong to $\GF^2$
and such that $Gf=g$ with $G \in (L^\infty)^{2 \times 2}$.
If $\GS:={\frak s}_{(\det G) f^T,-g^T}= \{0\}$,
then $\ker A^{\theta,\alpha}_\phi =  \GK f_1$, where
\beq\label{eq:7.9}
\GK =   \tilde f^T [ (H^2_+)^2]_{f \tilde f^T}
\cap
\tilde g^T [ (H^2_-)^2] _{g \tilde g^T}.
\eeq
\end{thm}

We now apply Theorem \ref{thm:7.1} and the results of the previous sections
to study the behaviour of an ATTO of the form \eqref{eq:7.5} with finite rank. 
As explained in \cite[Sec. 3]{CP17} we are led to take
the symbol to be 
\beq\label{eq:7.9bis}
\phi= \bar\theta R_+ - \alpha R_- +
\sum_{j=1}^N \frac
{\bar\theta P^\alpha_{n_j-1}(t_j)-\alpha P^{\bar\theta}_{n_j-1}(t_j)}
{(z-t_j)^{n_j}},
\eeq
where 

(i) $R_\pm$ are rational functions vanishing at $\infty$, such that $R_-$
has no poles in $\DD^e \cup \TT$ and $R_+$ has no poles in $\DD \cup \TT$;

(ii) $t_j \in \TT$ $(j=1,\ldots,N)$ are regular points for $\theta$ and $\alpha$, i.e., 
$\theta$ and $\alpha$ are analytic in a neighbourhood of each $t_j$ (in which case
$\bar\theta$, which can be extended to a  definition outside $\overline\DD$  by $\bar\theta(z)= \overline{\theta(1/\bar z)}$ for 
$|z|>1$, is also analytic in a neighbourhood of $t_j$);

(iii) $P^\alpha_{n_j-1}$ and $P^{\bar\theta}_{n_j-1}$ are the Taylor polynomials of order
$n_j-1$, relative to the point $t_j$, for $\alpha$ and for $\bar\theta$,
respectively.

Defining
\begin{eqnarray}
R_2^+ &=& R_+ + \sum_{j=1}^N \frac
{ P^\alpha_{n_j-1}(t_j)}
{(z-t_j)^{n_j}} = \frac{Q_2}{\GE D_{2+}}, \label{eq:7.11}
\\
R_1^-  &=&  R_- - \sum_{j=1}^N \frac
{ P^{\bar\theta}_{n_j-1}(t_j)}
{(z-t_j)^{n_j}} = \frac{Q_1}{\GE D_{1-}},
\end{eqnarray}
where 
\beq\label{eq:7.13}
\GE=\prod_{j=1}^N (z-t_j)^{n_j}, \quad \hbox{with} \quad \deg\GE= n_\TT := \sum_{j=1}^N n_j,
\eeq
and $Q_1, Q_2$ are polynomials, $D_{2+}$ is the denominator of $R_+$, with $n_-$
zeros (including multiplicities) in $\DD^e$, and $D_{1-}$ is the denominator of $R_-$,
with $n_+$ zeros in $\DD$, we can write, from \eqref{eq:7.9bis},
\beq\label{eq:7.14}
\phi= \bar\theta R_{2+} - \alpha R_{1-}.
\eeq
Then $G$ in \eqref{eq:7.6bis} takes the form
\beq\label{eq:7.15}
G= \left[
\begin{matrix}
\bar\theta & 0 \\
\bar\theta R_{2+}-\alpha R_{1-} & \alpha 
\end{matrix}
\right]
\eeq
and one can verify that $f,g \in \GF^2$, defined by
\begin{eqnarray} 
f = \left[
\begin{matrix}f_1 \\ f_2
\end{matrix}
\right] &=&
\left[
\begin{matrix}\theta \\ \theta R_{1-}
\end{matrix}
\right] =
\theta
\left[
\begin{matrix}1 \\ \frac{Q_1}{\GE D_{1-}}
\end{matrix}
\right],
\nonumber
\\
g =
\left[
\begin{matrix}g_1 \\ g_2
\end{matrix}
\right] &=& 
\left[
\begin{matrix} 1 \\   R_{2+}
\end{matrix}
\right] =
\left[
\begin{matrix}1 \\ \frac{Q_2}{\GE D_{2+}}
\end{matrix}
\right],
\label{eq:7.16}
\end{eqnarray}
satisfy $Gf=g$ and have
left inverses, respectively,
\beq
\tilde f^T=[\bar\theta \quad 0 ], \qquad \tilde g^T=[1 \quad 0].
\eeq

To apply Theorem \ref{thm:7.1},
we start by showing that the assumption $\GS=\{0\}$ in Theorem \ref{thm:7.1}
is satisfied when $K_\alpha$ is ``large enough''.

\begin{prop}
Let $G$ be given by \eqref{eq:7.15} with $\alpha$ such that 
$\dim K_\alpha \ge m:= n_+ + n_- + n_\TT $, 
and let $f,g$ be defined as in \eqref{eq:7.16}.
Then $\GS= {\frak s}_{(\det G) f^T,-g^T}=\{0\}$.
\end{prop}
\beginpf
We have $\det G=\alpha \bar\theta$ and $\GS$ is defined by
\begin{eqnarray}
&&\alpha \left( \phi_{1+}+ \frac{Q_1}{\GE D_{1-}}\phi_{2+}\right) = \phi_{1-}+\frac{Q_2}{\GE D_{2+}}\phi_{2-} \label{eq:7.18} \\
& \iff &
\frac{\alpha}{D_{1-}} 
\underbrace{(\GE D_{1-}\phi_{1+}+Q_1 \phi_{2+})}_{\psi_+} = \frac{z^{n_-+n_\TT}}{D_{2+}}
\underbrace{\frac{\GE D_{2+}\phi_{1-}+Q_2 \phi_{2-}}{z^{n_-+n_\TT}}}_{\psi_-}
\nonumber\\
&\iff& \alpha D_{2+}\psi_+ = z^{n_-+n_\TT}D_{1-}\psi_-,
\nonumber
\end{eqnarray}
where $\psi_\pm \in H^2_\pm$. So $\psi=\psi_++\psi_- \in {\frak s}_{\alpha D_{2+},-z^{n_-+n_\TT}D_{1-}}=\{0\}$,
by Corollary \ref{cor:2.3bis} and Proposition \ref{nprop:6.9} (i)
(with $p_1=D_{2+}$, $p_2=-z^{n_-+n_\TT}D_{1-}$,
$n_{2\DD}=n_++n_-+n_\TT$ and $ n_{1\TT}=n_{1\DD}=0$) and it follows that the left-hand side of \eqref{eq:7.18} is
\[
\alpha\left( \phi_{1+}+ \frac{Q_1}{\GE D_{1-}} \phi_{2+} \right) = \frac{\alpha}{\GE D_{1-}}\psi_+ = 0
\]
and, of course, analogously for the right-hand side of \eqref{eq:7.18}, so $\GS=\{0\}$.
\endpf

Next we characterise the spaces $[(H^2_+)^2]_{f\tilde f^T}$ and
$[(H^2_+)^2]_{g\tilde g^T}$ in \eqref{eq:7.9}. Note that
\beq
f \tilde f^T = \left[
\begin{matrix}
1 & 0 \\ R_{1-} & 0
\end{matrix}
\right]
\qquad \hbox{and} \quad
g\tilde g^T = \left[
\begin{matrix} 1 & 0 \\ R_{2+} & 0 
\end{matrix}\right],
\eeq
so
\beq\label{eq:7.19bis}
(\phi_{1+},\phi_{2+}) \in [(H^2_+)^2]_{f\tilde f^T} \iff
R_{1-}\phi_{1+} \in H^2_+ \iff \phi_{1+} \in \GE D_{1-}H^2_+
\eeq
and
\beq\label{eq:7.20bis}
(\phi_{1-},\phi_{2-}) \in [(H^2_-)^2]_{g\tilde g^T} \iff
R_{2+}\phi_{1-} \in H^2_- \iff
\phi_{1-} \in \frac{\GE D_{2+}}{z^{n_-+n_\TT}} H^2_-.
\eeq
We can now formulate the main result of this section.

\begin{thm}
Let $\phi$ be given by \eqref{eq:7.14} and \eqref{eq:7.11}--\eqref{eq:7.13},
and
let $\dim K_\alpha \ge m:=n_++n_-+n_\TT$.
Then $\ker T_G$ and $\ker A^{\theta,\alpha}_\phi$ do not depend on $\alpha$ and we have
\[
\ker T_G = \GE D_{1-}D_{2+} \ker T_{\bar\theta z^m} 
\left[ \begin{matrix} 1 \\ R_{1-} \end{matrix} \right] 
\quad \hbox{and} \quad
\ker A^{\theta,\alpha}_\phi = \GE D_{1-}D_{2+} \ker T_{\bar\theta z^m}.
\]
This holds, in particular, for any infinite-dimensional $K_\alpha$.
\end{thm}

\beginpf
By \eqref{eq:7.9}, \eqref{eq:7.19bis} and \eqref{eq:7.20bis}\begin{eqnarray*}
\GK &=& [\bar\theta \quad 0]
\left[ \begin{matrix} \GE D_{1-}H^2_+\\  H^2_+ \end{matrix} \right]
\cap [1 \quad 0]
\left[ \begin{matrix}  \frac{\GE D_{2+}}{z^{n_- + n_\TT}}H^2_-\\  H^2_- \end{matrix} \right] \\
&=& \bar\theta \GE D_{1-}H^2_+ \cap \frac{\GE D_{2+}}{z^{n_- + n_\TT}}H^2_-\\
&=&\bar\theta\GE D_{1-} \k^+_{\bar\theta D_{1-},-\frac{  D_{2+}}{z^{n_- + n_\TT}}}\\
&=& \bar\theta\GE D_{1-} \k^+_{\bar\theta D_{1-}z^{n_- + n_\TT},-   D_{2+}}\\
&=& \bar\theta\GE D_{1-}D_{2+} \k^+_{\bar\theta z^{n_+ + n_- + n_\TT},-   1}\\
&=& \bar\theta \GE D_{1-}D_{2+} \ker T_{\bar \theta z^m},
\end{eqnarray*}
where we used Proposition \ref{prop:6.12bis}, and the result follows from Theorem
\ref{thm:7.1}.
\endpf

Note that, if $z^m$ divides $\theta$ then $\ker T_{\bar\theta z^m}$ is the
model space $K_{\theta/z^m}$, so $\ker A_\phi^{\theta,\alpha}$ is isomorphic to
$K_{\theta/z^m}$.

\section*{Acknowledgements}
The authors are grateful to Yuxia Liang and an anonymous referee for
helpful comments.

\end{document}